\documentclass[12pt]{amsart}
\usepackage{amsmath,amssymb,amsthm,amscd,amsopn,amsfonts,amsxtra} 
\usepackage{latexsym,array,color}
\usepackage[mathscr]{eucal}
\usepackage{graphicx,psfrag}
\usepackage{epsfig}
\usepackage{color}
\newtheorem{theorem}[subsection]{Theorem} 
\newtheorem{lemma}[subsection]{Lemma}

\newtheorem{proposition}[subsection]{Proposition}
\newtheorem{assumption}[subsection]{Assumption}

\newtheorem{remark}[subsection]{Remark}

\definecolor{orange}{rgb}{0.995, 0.75, 0.35}
\definecolor{purple}{rgb}{0.7, 0.2, 0.5}
\definecolor{royalblue}{rgb}{0.2, 0.7, 0.8}
\def\al{\alpha}

\def\eps{\epsilon}

\def\lam{\lambda}

\def\vphi{\varphi}
\def\De{\Delta}

\def\sh{\sinh}

\def\sech{\mathrm{sech}}
\def\supp{\mathrm{supp}}

\newcommand{\cal}{\mathcal}

\newcommand{\la}{\langle}
\newcommand{\ra}{\rangle}
\newcommand{\nd}{\noindent}

\newcommand{\vs}{\vspace}

\newcommand{\hB}{\hfill$\Box$}

\newcommand{\Z}{\mathbb{Z}}
\newcommand{\R}{\mathbb{R}}

\newcommand{\N}{\mathbb{N}}

\begin{document}
\title[Littlewood-Paley theorem for Schr\"odinger operators]
{Littlewood-Paley theorem for Schr\"odinger operators}
\author{Shijun Zheng}
\address{Department of Mathematics \\
               University of South Carolina\\
               Columbia, SC 29208}
\email{shijun@math.sc.edu}
 \urladdr{http://www.math.sc.edu/~shijun}
\thanks{This work is supported by DARPA grant HM1582-05-2-0001. The author gratefully
thanks the hospitality and support of Department of Mathematics, University of
South Carolina, during his visiting at the Industrial Mathematics Institute} 

 %35Q40}
\date{July 26, 2006}

%\tableofcontents 
%\input{abstr}
\begin{abstract}

Let $H$ be a Schr\"odinger operator on $\R^n$.
%Let  $H$ denote a selfadjoint operator whose heat kernel and its derivative 
%satisfy the upper Gaussian bound.  
Under a polynomial decay condition for the kernel of its spectral operator,   
%certain decay,  
 we show that the %homogeneous 
Besov spaces and Triebel-Lizorkin spaces  
associated with $H$ are well defined.  
We further give a Littlewood-Paley characterization of $L_p$ spaces 
as well as Sobolev spaces in terms of dyadic functions of $H$. 
%a selfajoint operator %where $\mL$ is a elliptic operator, 
% Schrodinger operator with magnetic potential.
This generalizes and strengthens the previous result when the heat kernel of $H$ 
satisfies certain upper Gaussian bound. %Hermite operator 
%\cite{E95,D97} 
%as well as other Schr\"odinger operators with short-range potentials.

%especially for negative potentials. main approach and idea 
%extend to general potentials.  Hermite  Laguerre 
\end{abstract}

\maketitle 

\section{Introduction and main results}\label{S1}

Recently the theory of function spaces associated with Schr\"odinger operators have been drawing attention
in the area of harmonic analysis and PDEs \cite{JN94,DP05,BZ05,OZ,Sch05b,E95,D97,DZ02,DZ05,DY05}.
In \cite{E95,D97,BZ05,OZ}  it is proved that the Besov and Trieble-Lizorkin spaces associated with a Schr\"odinger operator are well defined,  in some particular cases. 
  In this note we aim to extend the result for general Schr\"odinger operators on $\R^n$. 
Furthermore we are interested in obtaining a  Littlewood-Paley decomposition for 
the $L_p$ spaces as well as Sobolev spaces using dyadic functions of  $H$.  

Let $H=-\De+V $ be a Schr\"{o}dinger operator 
that is selfadjoint in $L_2(\R^n)$ with
 a real-valued potential function $V$. 
Then for a Borel measurable function $\phi$,   one can define the 
spectral operator $\phi(H)$ by functional calculus
$\phi(H)=\int_{-\infty}^\infty \phi(\lam)dE_\lam   $, where $dE_\lam$ is 
the spectral measure of $H$.  The kernel of $\phi(H)$ is denoted $\phi(H)(x,y)$.
%${\cal D}(\phi(H))=  \{u\in {\cal H}:\int |\phi(\lam)|^2 d(E(\lam)u,u)<\infty\}$.

Let $ \{\varphi_j\}_{j\in\Z}\subset C_0^\infty ({\R}) $ be a smooth dyadic system satisfying the conditions 
%\mbox{$\forall x \in \mathbf{R}$,}
\quad $\text{(i)}\; %\supp \;\Phi \subset \{ x: |x|\le 1\},
 \supp\; \varphi_j
\subset \{ x: 2^{j-2}\le |x|\le 2^j\} $ \quad 
 
  %\hspace{.12in} 
 $\text{(ii)}\; |\vphi_j^{(k)}(x)|\le c_k 2^{-kj}\, ,    \qquad \forall j\in \Z$, $k\in \N_0=\{0\}\cup\N$
%\; \textrm{there is}\; c >0 \;  \textrm{such that}\; \vert \Phi (x)\vert
%\ge c \;\textrm{if} \;|x| \le 5/6 \;\textrm{and}\;  \vert \varphi (x)\vert
%\ge c \;\; \textrm{if} \; \\  3/10 \le |x| \le 5/6; $
 
  (iii)  \begin{align*}
\sum_{j=-\infty}^\infty \varphi_j(x) \approx c>0\, ,  \;  \qquad \forall x\neq 0\;.     \label{eq:dyadic-id}
%&\Phi(x)  +\sum_{j=1}^\infty \varphi_j(x) \approx c>0 
\end{align*}

Let $0< p<\infty$, $0< q\le\infty$ and $\al\in\R$. The  homogenous {\em Triebel-Lizorkin space} $\dot{F}^{\al,q}_p(H)$  
%associated with $H$ 
is defined as the completion of the Schwartz class $\cal{S}(\R^n)$  %linear space 
%$S_H:=\{ \chi(H)u:  u\in C^\infty(\R^n), \chi\in C^\infty_0(\R) \}$, 
with the quasi-norm   % $\Vert f\Vert_{F^{\al,q}_p(H)}$ is given by
$$\Vert f\Vert_{\dot{F}^{\al,q}_p(H)} = 
 \Vert\big(\sum_{j=-\infty}^\infty 2^{j\al q} \vert \vphi_j(H)f(\cdot)\vert^q\big)^{1/q}\Vert_p.
$$

\nd Similarly, if $0<p\le\infty$, $0< q\le\infty$, the homogeneous {\em Besov space} $\dot{B}^{\al,q}_p(H)$ 
%associated with $H$, denoted by   , 
is defined by the quasi-norm
$$\Vert f\Vert_{\dot{B}^{\al,q}_p(H)} = 
 \big(\sum_{j=-\infty}^\infty 2^{j\al q} \Vert \vphi_j(H)f \Vert_p^q\big)^{1/q}.$$

%\rk the subsp $S_H$ can be chosen to be S(\R^n) or C^\infty_0,
% which imply, since C^\infty_0\subset S(\R^n)\subset B^s(H), some s>0  
%we have
%\[   (B^s(H) )'\subset S'\subset D'  \]
%So it is appropriate to choose S_H= C^\infty_0 or S

%\begin{CD}
%B^{\al,q}_p(H)@> = >> (L^p, \mathcal{D}(H^m) )_{\theta, q} \\
%@AAA         @AAA\\  
 %B^{2\al,q}_p(\R^d)@> =>> (L^p, W_p^{2m} )_{\theta, q} \\
%\end{CD}
%$$  
%\vspace{.1770in}
%where $ \theta =\frac{\al}{m}$. 

Throughout this note we assume $H$ satisfies the following: %assumption %\ref{a:phi-dec} 
%on the decay for kernel of $\phi(H)$ when $\{\phi_j\}_j$ is a dyadic system as above. 

\begin{assumption}\label{a:phi-dec}  
Let $\phi_j\in C_0^\infty(\R)$ be as in condition (i), (ii).  %$\supp\; \phi\subset [-1,1]$. 
Then  %$\phi_j(H)$
 for every $N\in\N_0$ there exists a constant 
 $c_N>0$ such that for all $j\in\Z$
\begin{equation} \label{e:ker-phi}
\vert  \phi_j(H)(x,y)\vert 
\le c_N\frac{ 2^{nj/2}}{(1+2^{j/2}|x-y|)^N}
\end{equation}
\begin{equation}\label{e:der-ker-phi}
\vert  \nabla_x \phi_j(H)(x,y)\vert 
\le c_N\frac{2^{(n+1)j/2}}{(1+2^{j/2}|x-y|)^N}.
\end{equation}
\end{assumption}

This is the case when $H$ is the Hermite operator $-\De+|x|^2$, 
or more generally, whenever $V$ is nonnegative and 
$H$  satisfies the upper Gaussian bound for the heat kernel and its derivative
 (see Proposition \ref{pr:etH-phi}).
  However, %in dimensions  one and two 
when the potential $V$ is negative, such a heat kernel estimate is {\em not}
available. Therefore it is necessary to consider  a more general condition as given 
in Assumption \ref{a:phi-dec}. 

Define the {Peetre maximal function} for $H$ as:  for $j\in\Z$, $s>0$

$$
\vphi_{j,s}^*f(x) = \sup_{t\in \R^n} \frac{| \vphi_j(H)f(t)|}{
 ( 1+2^{j/2}|x- t  |)^{s}} \, , \;\;%s>0
 $$

and 

$$
\vphi_{j,s}^{**} f (x) = \sup_{t\in \R^n} \frac{| (\nabla_t\vphi_j(H)f)(t)|}{
 ( 1+2^{j/2}|x -t  |)^{s}} \, .$$

The following theorem is a maximal characterization of the {\em homogeneous} 
%$F$ and $B$ 
spaces. 
By $\Vert \cdot\Vert_A \approx \Vert \cdot \Vert_B $ we mean equivalent norms. 

\begin{theorem}\label{th:homog-B-F} Suppose $H$ satisfies Assumption \ref{a:phi-dec}.
%Let $0< p\le\infty$, $0< q\le \infty$, $\al\in \R$. 
%Then 

a) If $0< p\le\infty$, $0< q\le \infty$, $\al\in \R$ and  $s>n/ p$, then
\[ \Vert f\Vert_{\dot{B}_p^{\alpha,q}(H)}\approx
\Vert \{2^{j\al}\vphi^*_{j,s}(H)f\}\Vert_{\ell^q(L_p)} \;. \]

b) If  $0< p< \infty$, $0< q\le \infty$, $\al\in \R$ and $s>n/\min (p,q)$, then 
\[
 \Vert f\Vert_{\dot{F}_p^{\alpha,q}(H)}\approx 
\Vert \{2^{j\al}\vphi^*_{j,s}(H)f\}\Vert_{L_p(\ell^q)} \;.
\]
\end{theorem}

It is well-known that such a characterization implies that any two dyadic systems satisfying (i), (ii), (iii)
 give rise to equivalent norms on $\dot{F}_p^{\alpha,q}(H)$ and $\dot{B}_p^{\alpha,q}(H)$. 
The analogous result also holds for the inhomogenous spaces $F^{\al,q}_p(H)$, $B^{\al,q}_p(H)$.
However, the homogeneous spaces,  %$L_p$ and Besov spaces, 
which cover both high and low energy portion of $H$, 
are essential and more useful  in proving Strichartz inequality for wave equations \cite{Sch05b,KT98}.
This is one reason of our motivation.  

Following the same idea in \cite{OZ}, using Calder\'on-Zygmund decomposition and Assumption
\ref{a:phi-dec}  we show that  $L_p(\R^n)=F^{0,2}_p(H)$ 
 if $1<p<\infty$.  We thus obtain the Littlewood-Paley theorem for $L_p$ spaces. 
\begin{theorem}\label{th:L-P-Lp}   Suppose $H$ satisfies Assumption \ref{a:phi-dec}.
 If $1<p<\infty$, then
\begin{equation*}
\Vert f\Vert_{L_p(\R^n)} \approx %\Vert f\Vert_{\dot{F}_p^{\alpha,q}}\sim 
\Vert \big(\sum^\infty_{j=-\infty} |\vphi_j(H)f(\cdot)|^2\big)^{1/2} \Vert_{L_p(\R^n)}\;.
\end{equation*}
\end{theorem}
%\rk
 Under additional condtion on $V$, e.g., $|\partial_x^kV(x)|\le c_k$, $|k|\le 2m_0-2$ 
for some $m_0\in \N$, %(or $\partial^k_x V$ is a space multiplier)
we can %using lifting properties of $H$ to 
characterize the Sobolev spaces $H^{2s}_p(\R^n)=F^{s,2}_p(H)$, $1<p<\infty$, $| s|\le m_0$
%$s\in\R$ 
with equivalent norms
\begin{equation*}
 \Vert f\Vert_{H_p^{2s}(\R^n)} \approx %\Vert f\Vert_{\dot{F}_p^{\alpha,q}}\sim 
\Vert \big(\sum^\infty_{j=-\infty} 2^{2js}|\vphi_j(H)f(\cdot)|^2\big)^{1/2} \Vert_{L_p(\R^n)}\;.
\end{equation*}

\section{Proofs of Theorem \ref{th:homog-B-F} and Theorem \ref{th:L-P-Lp}}

The proof of Theorem \ref{th:homog-B-F} is standard and follows from 
 Bernstein  type inequality (Lemma \ref{l:bern-der-phi}) and
  Peetre type maximal inequality (Lemma \ref{l:phi*-M})  for maximal functions. 
% Both lemmas rely on %and essentially depends on 
%the estimates (\ref{e:ker-phi}), (\ref{e:der-ker-phi})   for $\vphi_j(H)$ 

\begin{lemma}\label{l:bern-der-phi} %(Bernstein) 
For $s>0$, there exists a constant $c_{n,s}>0$ such  that
 for all $j\in\Z$ 
\[
\vphi_{j,s}^{**}f(x) \leq c_{n,s} 2^{j/2} %\max_{\substack{0\leq k \leq 2N\\* \pm} } 
\vphi_{j,s}^{*} f( x), \qquad \quad \; \forall f\in \cal{S}(\R^n).
\]
\end{lemma}

Similar to \cite{Tr83,OZ}  Lemma \ref{l:bern-der-phi} can be easily proved using (\ref{e:der-ker-phi}) with %$\al=1$ 
$N> n+s$ 
and the identity 
\[  \vphi_j(H)f(x) =%\sum_{\nu = -1}^1  (\phi \psi)_{j+\nu}
\psi_j(H) \vphi_j(H)f , \]
%\;\; \qquad f\in \cal{S}(\R^n),   $$
where $\psi_j(x)=\psi(2^{-j}x)$ with $\psi\in C^\infty_0$, $\supp\;\psi\subset \{ \frac{1}{5}\le |x|\le \frac{5}{4} \}$
and $\psi(x)=1$ on $\{ \frac{1}{4}\le |x|\le 1 \}$.  

Let $M$ denote the Hardy-Littlewood maximal function
\begin{equation}\label{e:F-S-max}
 M f (x)=  \sup_{B\ni x}\frac{1}{|B|} \int_{B}  |f(y)| \, d y
\end{equation}
where the supreme is taken over all %geodesic
 balls $B$ in $\R^n$ centered at $x$. 

\begin{lemma}\label{l:phi*-M}   %(Peetre's maximal inequality)
Let $0<r<\infty$ and $s=n/r$. 
 %If $f\in \cal{S}(\R^n)$, 
 Then 
%\begin{itemize}
%\item[(a)] 
for all $j\in\Z$
\begin{equation}\label{e:phi*-M} 
 \vphi_{j,s}^*f(x) \le c_{n,r} [M({|\vphi_{j} (H)f|}^r)]^{1/r}(x),\;\qquad \; \forall f\in \cal{S}(\R^n) .
 %0<r<\infty
\end{equation}
\end{lemma}

\begin{proof}  %First note that since $\partial_x e^{-t_jH}(x,y)$ exist, 
 Let $g(x) \in C^1(\R^n)$. As in \cite{Tr83,BZ05}, the mean value theorem
gives for $z_0\in \R^n$,  $\delta>0$ 

$$
|g(z_0)|\leq \delta \sup_{\substack{|z-z_0|\leq \delta}} |\nabla g(z)|
 +c_{n,r}\delta^{-n/r}\left( \int_{|z-z_0|\leq \delta} |g|^r dz \right)^{1/r}.$$

\nd
Put $g(z)=\vphi_j(H)f( x -z)$ %\in C^1$ 
to get  %(with $0<\delta\leq 1)$,   

\begin{align*}
&\frac{|\vphi_j(H) f ( x-z)|}{(1+2^{j/2}|z|)^{n/r}}
\leq \delta \sup_{\substack{|u-z|\leq \delta}} \frac{(1+2^{j/2}|u|)^{n/r}
|\nabla (\vphi_j(H)f)(x-u)|}{(1+2^{j/2}|z|)^{n/r}
(1+2^{j/2}|u|)^{n/r}} \\
+& c_{n,r}\delta^{-n/r} (1+ 2^{j/2}|z|)^{-n/r} (\int_{|u-z|\leq \delta} 
|\vphi_j(H)f(x  -u)|^r du)^{1/r}\\ 
%\leq& \delta (1+2^{j/2}\delta)^{n/r} \sup_{\substack{u\in \R^n}}
%\frac{ | \nabla(\vphi_j(H)f)( x  -u)|}{(1+2^{j/2}
%|u|)^{n/r}}   \; {\color{red}checked} \\
 %+& 
%c\delta^{-n/r} (1+2^{j/2}|z|)^{-n/r} (\int_{|u|\leq |z|+\delta} 
%|\vphi_j(H)f(x  -u)|^r du)^{1/r}\\
\le& \delta (1+2^{j/2}\delta)^{n/r} \vphi_{j,s}^{**}f(x) +c_{n,r}{\delta}^{-n/r} 
\big{(}\frac{|z|+\delta}{1+2^{j/2}|z|}\big{)}^{n/r} [M(|\vphi_j(H)f |^r ) (x)]^{1/r} \\  %\;\; {\color{red}checked}
 \le& c_{n,r}\delta (1+2^{j/2}\delta)^{n/r} 2^{j/2} \vphi_{j,s}^{*}f(x) +c_{n,r}{\delta}^{-n/r} 
\big{(}\frac{|z|+\delta}{1+2^{j/2}|z|}\big{)}^{n/r} [M(|\vphi_{j}(H)f |^r )]^{1/r}(x)\\
\le&c_{n,r} \epsilon(1+\eps)^{n/r} \vphi_{j,s}^*f(x) +c_{n,r} (1+ \epsilon^{-1})^{n/r}
 [M( |\vphi_{j}(H)f |^r )]^{1/r} (x),
%(0<\epsilon \leq 1)
\end{align*}
by setting $\delta= 2^{-j/2}\epsilon $, $\eps>0$ and using Lemma \ref{l:bern-der-phi}. 
Finally, taking $\eps>0$ sufficiently small 
establishes (\ref{e:phi*-M}).  
\end{proof}
 
  Now Theorem \ref{th:homog-B-F} is a consequence of
  Lemma \ref{l:phi*-M} and  the following well-known lemma %\ref{l:F-S-max} 
on Hardy-Littlewood maximal function  by a standard argument; 
see \cite{Tr83} or \cite{E95,OZ} for some simple details. 

%The following is a vector valued inequality for H-L maximal function on $\R^n$. %\cite{Skryz}. 

\begin{lemma}\label{l:F-S-max}  a) If $1<p\le\infty$, then
\begin{equation}
  \Vert  M f\Vert_{L_p(\R^n)}\le C_p\Vert f\Vert_{L_p(\R^n)} \;.
\end{equation}
b) If $1<p<\infty$, $1<q\le\infty$, then 
\begin{equation}
\Vert \bigg(\sum_j | M f_{j}|^q\bigg)^{1/q}\Vert_{L_p(\R^n)}
\le C_{p,q}\Vert \bigg(\sum_j |f_{j}|^q \bigg)^{1/q}\Vert_{L_p(\R^n)} \;.\label{e:F-S-max}
\end{equation}
\end{lemma}

%{\bf Conclusion:} If low energy becomes a problem, it  seems 
%to suggest that we need to modify the definition of maximal function $\phi_j^*f$.

\subsection{Proof of Theorem \ref{th:L-P-Lp} }  From the proof of the identification of $F^{0,2}_p(H)$ spaces
\cite[Theorem 5.1]{OZ} we observe that the estimates 
in (\ref{e:ker-phi}), (\ref{e:der-ker-phi}) imply %that    
\begin{equation}\label{e:f-F-Lp}
\Vert f\Vert_{F^{0,2}_p(H)} \approx   \Vert f\Vert_{L_p} \;, \qquad \; 1<p<\infty
\end{equation}
by applying $L_p(\ell^2)$-valued Calder\'on-Zygmund decomposition. 
On the other hand,  Theorem \ref{th:homog-B-F} suggests %in particular with $\al=0$ 
that
\begin{equation}\label{e:f-F-Lplq}
  \Vert f\Vert_{F^{\al,q}_p(H)} \approx \Vert \{2^{j\al}\vphi_{j}(H)f\} \Vert_{L_p(\ell^q)}  
\end{equation}
 whenever $\{\vphi_j\}_{j\in\Z}$ is a dyadic system satisfying (i), (ii), (iii).  

Combining (\ref{e:f-F-Lp}) and (\ref{e:f-F-Lplq}) with $\al=0$, $q=2$ proves Theorem \ref{th:L-P-Lp}.
\hB

  %Definition. A function $f$  is said to be an element of $H_A^1$ if the
%maximal function ${\cal M}f(x):= \sup_{t>0} |e^{-tA}f(x) |$ belongs to $L^1$. 
%Under the above condition on $V$.  

\vs{.172in}
\nd
\begin{remark}   For $p=1$,  Dziuba\'nski and  Zienkiewicz \cite{DZ05} recently obtained a characterization of Hardy space associated with
$H$ 
and showed that if a compactly supported positive potential $V$ is in $L^{n/2+\eps}$, 
$n\ge 3$, then  %and is compacted supported,  
\[
\left\|f\right\|_{\cal{H}^1}\approx  \left\|wf\right\|_{H^1\left({\Bbb R}^d\right)} \;,
\] 
where $\cal{H}^1=\{f\in L^1: \sup_{t>0}  |e^{-tH}f(\cdot) |  \in L^1  \}$ and the weight $w$ is defined by $w(x)= \lim_{t\to\infty} \int_{\R^n} e^{-tH}(x,y)  dy$. 
%$\|g\|_{H^1({\Bbb R}^d)}= \|\sup_{t>0}|P_t g(x)|\|_{L^1}$ is the norm in the classical real Hardy  H^1({\Bbb R}^d)         
 It  would be very interesting to see whether one can give a Littlewood-Paley characterization 
of $\cal{H}^1$ in the sense of Theorem \ref{th:L-P-Lp}. 
\end{remark}
									 	
\section{Potentials satisfying upper Gaussian bound}

In this section we show that %give a condition on nonnegative potentials V so that 
 Assumption \ref{a:phi-dec} is verified when $H$ satisfies the upper Gaussian bound 
 (\ref{e:der-etH-gb}) for its heat kernel.  We begin with a weighted $L^1$ inequality,  %for $g(2^{-j}H)$,  
which is an easy consequence of
\cite[Lemma 8]{He90a} by a %regular
 scaling argument. 
%Observe from \cite[Lemma 9]{He90a} that a regular scaling argument gives
\begin{lemma} (Hebisch)\label{l:sca-heb-dec}  
Suppose $V\ge 0$ and 
 $e^{-tH}$ satisfies 
%the upper Gaussian bound 
\begin{equation}\label{e:etH-gb}
0\le e^{-tH}(x,y)\le c_n t^{-n/2} e^{-c |x-y|^2/t}\;,\qquad \forall t>0.
\end{equation}
%where $C=C(V)$ d.n. depend the scaling V\to \al^2 V(\al x)
%e.g. V positive or $\Vert V_-\Vert<c$ 
 If  $s> (n+1)/2+\beta$,  $\beta\ge 0$ and %$g\in H^s(\R)$
    $\supp\; g\subset [-10,10]$, then   %$j\in\Z$. 
 \begin{align*}
\sup_{j\in\Z, \, y\in\R^n}\Vert g(2^{-j}H) (\cdot,y) \la 2^{j/2} (\cdot-y)\ra^\beta \Vert_{L^1(\R^n)} 
\le c_n \Vert g\Vert_{H^s(\R)} \, ,%2^{j/2} (1+ 2^{j/2}|x-y|)^{-N} .
\end{align*}
 where $\la x\ra:=1+|x|$ and $\Vert\cdot\Vert_{H^s}$ denotes the 
usual Sobolev norm.
\end{lemma} 

\begin{remark} It is known that (\ref{e:etH-gb}) holds whenever $V\ge 0$ is 
locally integrable. %in $L^1_{loc}(\R^n)$  \cite{Ou06}
\end{remark}

\begin{proposition}\label{pr:etH-phi}   Let $\al=0$, $1$. 
Suppose $V\ge 0$ and $e^{-tH}$ satisfies the upper Gaussian bound %for $\al=0,1$
\begin{equation}\label{e:der-etH-gb}
|\nabla^\al_x e^{-tH}(x,y)| \le c_n t^{-(n+\al )/2} e^{-c  |x-y|^2/t}\,, \qquad \forall t>0.
\end{equation}
 If $\{\vphi_j\}_{j\in\Z}$ is  %admissible  
 a dyadic system  
 satisfying (i), (ii), %(x)=\phi(2^{-j}x)$, $\phi\in C^\infty_0(\R)$, $j\in\Z$. 
 then for each $N\ge 0$ 
\begin{align*}
\vert \nabla_x^\al \vphi_j(H) (x,y) \vert 
\le c_N 2^{j(n+\al )/2} (1+2^{j/2}|x-y|)^{-N} , \qquad \forall j\,.
\end{align*}
\end{proposition} 

\begin{proof}  Write
\begin{align*}
&\nabla^\al_x \vphi_j(H)(x,y)=\int_z \nabla_x^\al e^{-tH}(x,z) (e^{tH}\vphi_j(H) )(z,y) dz.
\end{align*}
%Setting $t_j = 2^{-j}$ and using 
By (\ref{e:der-etH-gb}) we have 
\begin{align*}
&|\nabla^\al_x \vphi_j(H)(x,y)|\\
%\lesssim& \int_z  t_j^{-(n+\al )/2} e^{-c|x-z|^2/t_j} \la x-z\ra^N \la x-z\ra^{-N} |(e^{t_jH}\vphi_j(H) )(z,y)| dz\\ 
\le c_n&  t^{-(n+\al )/2} \int e^{-c|x-z|^2/t} \la(x-z)/\sqrt{t}\ra^N 
\la(x-z)/\sqrt{t}\ra^{-N} \la(z-y)/\sqrt{t}\ra^{-N} \\
&\qquad \qquad \quad \cdot \la(z-y)/\sqrt{t}\ra^{N} |(e^{tH}\vphi_j(H) )(z,y) | dz\\ 
\le c_n& t^{-(n+\al )/2} \la(x-y)/\sqrt{t}\ra^{-N} \int \la(z-y)/\sqrt{t}\ra^{N} | (e^{tH}\vphi_j(H) )(z,y)| dz.
%\le& c_{N,n}t_j^{-(n+\al )/2}  \la(x-y)/\sqrt{t_j}\ra^{-N} 
%\Vert \la(\cdot -y)/\sqrt{t_j}\ra^{N+n/2+\eps} (e^{t_jH}\vphi_j(H) )(\cdot,y) \Vert_{L^1(\R^n)}
\end{align*}
%According to a weighted ineq  \cite[Lemma 9]{He90a}  %[DOS02]
 %\cite[(12)]{He90a}
 %\[ g(2^{-j}H)(x,y)\le c \Vert g\Vert_{H^N} (1+|x-y|)^m  \qquad N=n/2+ m+1 \]

Setting $t=t_j:=2^{-j}$,  we see that   %it is easy to verify 
$g_j(x):=e^{t_jx}\vphi_j(x)  $  also satisfies conditions (i), (ii).  
%Let $g_0(x)=g_j(2^j x)$. 
Writing  $g_j(x)=g_0(2^{-j}x)$, then $\supp\; g_0\subset \{\frac{1}{4}\le |x|\le1\}$ and  
%$g_j(H)(x,y)\le c\Vert g_0\Vert_{H^s} 2^{j/2} (1+ 2^{j/2}|x-y|)^{-N} $
%c_n 2^{j/2} (1+ 2^{j/2}|x-y|)^{-N} $.
\[ \Vert g_0\Vert_{H^N(\R)}\le \Vert g_j(2^j x)\Vert_{C^N(\R)} \le c_N \;.
\]
Thus an application of  Lemma \ref{l:sca-heb-dec} with $g=g_0$, $\beta=N$ proves the
proposition.
\end{proof}

\subsection{Hermite operator $H=-\De+|x|^2$}

To verifies Assumption \ref{a:phi-dec}   it is sufficient to show  
$e^{-tH}$ satisfies the upper Gaussian bound in (\ref{e:der-etH-gb}),   
according to Proposition \ref{pr:etH-phi}.

For $k\in\N_0$, let $h_k$ be the $k^{th}$ Hermite function with $\Vert h_k\Vert_{L_2(\R)}=1$
such that 
\[  (-\frac{d^2}{dx^2}+x^2) h_k = (2k+1)h_k\,.
\]
Then $\{h_k\}_0^\infty$ forms a complete orthonormal system (ONS) in $L_2(\R)$. 
%Hermite polynomials $H_0(x)=1$, $H_1(x)=2x$ ($H_{-1}(x)=0$). 
 % Hermite functions $h_0(x)=\pi^{-1/4}e^{-x^2/2}$, $h_1(x)=\pi^{-1/4}2^{1/2} xe^{-x^2/2}$
%($h_{-1}(x)=0$).  
In $L_2(\R^n)$, the ONS is given by $\Phi_k(x):= 
h_{k_1}\otimes \cdots \otimes h_{k_n}$, $k=(k_1,\dots,k_n)\in \N_0^n$. 

By Mehler's formula \cite[Ch.4]{Th93} or \cite{Th06},  the heat kernel %of $e^{-tH}$ 
has the expression 
\begin{align*}
e^{-tH}(x,y)=& \sum_{k\in\N_0^n} e^{-t (n+2|k|)} \Phi_k(x) \Phi_k(y)\\
=& \frac{1}{(2\pi \sh (2t) )^{n/2} } e^{-\frac{1}{2}\coth (2t) (|x|^2+|y|^2) + \mathrm{cosech (2t)} x\cdot y }
\end{align*}
for all $t>0$, $x, y\in \R^n$.

It is easy to calculate to find that there exist constants $c,c'>0$, 
$0<c_0,c_1,c'_0,c'_1<1$ and $t_0>1$ such that 
\[
p_{\frac{t}{2} }(x,y) \le c
\begin{cases} 
t^{-n/2} e^{-c_0 |x-y|^2/t}&\quad t\le t_0\\
e^{- n t /2} e^{-c_1 |x-y|^2} &\quad  t>t_0
\end{cases} \]
\[
|\nabla_x p_{\frac{t}{2} }(x,y) | \le c'
\begin{cases} 
t^{-(n+1)/2} e^{-c'_0 |x-y|^2/t} &\quad t\le t_0\\
e^{- n t /2} e^{-c'_1 |x-y|^2}  &\quad t>t_0\,,
\end{cases}
\]
where $p_t(x,y):=e^{-tH}(x,y)$. Hence (\ref{e:der-etH-gb}) holds. 
%noting $\ch t= 1+\frac{t^2}{2}+O(t^4)$, $t\to 0$, 
%$ \frac{\ch t}{\sh t} = 1+(2+0)e^{-2t}$, $t\to\infty$.

\begin{remark}  For the Hermite operator, %$H$, 
%Epperson \cite{E97} obtained 
the decay estimates similar to (\ref{e:ker-phi}), (\ref{e:der-ker-phi})  were previously obtained
in \cite{E97} in one dimension and \cite{D97} in $n$-dimension.
%of $\vphi_j(H)(x,y)$ for $\vphi_j\in C^\infty_0$ 
%in one dimension. Later his result was extended to the $n$-dimension in 
The latter used Heisenberg group method. 
 \mbox{Proposition \ref{pr:etH-phi} } shows that 
using heat kernel estimate we can obtain a simpler proof.   
\end{remark}
%Laguerre 

\begin{remark}
 When $V$ is negative,  
 the heat kernel estimate (\ref{e:etH-gb}) is not available, espcecially in low dimensions $n=1,2$. 
but Assumption \ref{a:phi-dec} still holds in the high energy case ($j\ge 0$) 
for certain short range potentials. %admitting only   
A special example is the one dimensional P\"oschl-Teller model 
$V(x)= -\nu(\nu+1)\,\sech^2 x$, $\nu\in \N$, cf.  \cite{OZ}.  
We will discuss the problem in more detail in %the paper 
 \cite{Z06} where 
$V$ is assumed to have only polynomial decay at infinity. 
\end{remark}

%cases where $C^\infty_0\subset B^s(-\De+V)$}

%A result on perturbation of order zero is also valid, that is  

%Prop. Let $V=V_+-V_-$, $V_+\in L^1_{loc}(\Om)$, $V_-\in L^1(\Om)$. 
%Then the above ineq holds, in particular, 
%\[  |e^{-tH}(x,y)| \le A t^{-s} \exp(-c (d_M(x,y))^{2m/(2m-1)} t^{-1/(2m-1)} + c' t), \]
%where $s=n/2m<1$. 

%\rk Thus if $n=1, m=1$  then $-d^2/dx^2+V$ satisfies the above ineq.
%which means if we want the usual Gaussian bound in one dimension
%($c'=0$) we have to make V less singular.  e.g. short-range.   or Kato class.

%[Barbatis, G.]   Sharp heat-kernel estimates for higher-order operators with singular coefficients 2004. 

%{\bf Acknowledgment}. \quad The author would like to thank {\mbox Grillakis,} Machedon, 
%Strichartz and Tao for numerous discussions and valuable comments.

\end{document}